%% file: uld-lattices.tex
\documentclass[11pt]{article}

\usepackage{a4wide}
\usepackage{amssymb}
\usepackage{amsmath}
\usepackage{amsthm}
\usepackage{url}

\usepackage{graphicx}
\usepackage{color}
\graphicspath{{Figures/}}
\input{psfig.sty}
\usepackage{wrapfig}

\usepackage{psfrag}

\def\comment#1{}

\def\NI{\noindent}
\def\ni{\noindent}

\def\sk{\smallskip}

\def\term#1{{\bf #1}\marginpar{\raggedright{\small\sf #1}}}
\def\ncterm#1{{\bf #1}} 
\def\nct#1{{\bf #1}} 

\def\ITEMMACRO #1 ??? #2 ???{\par\vskip4pt\noindent%
\hangindent=#2em\setbox0\hbox{#1\kern4pt}%
\ifdim\wd0<\hangindent\setbox0\hbox to\hangindent{\hss#1\kern7pt}\fi%
\box0\ignorespaces}

\def\Item(#1){\ITEMMACRO {\rm (#1)} ??? 1.8 ???}

\let\Bitem=\bItem
\def\BrackItem[#1]{\ITEMMACRO [#1] ??? 1.8 ???}

\def\Note{\ITEMMACRO {\NI\bf Note.}  ??? 1.8 ??? \\ \bgroup\sf}
\def\EndNote{\par\egroup\medskip\ni}

\newtheoremstyle{meiner} 
    {4pt}{3pt}           
    {\sffamily\upshape}  
    {}                   
    {\sffamily\bfseries} 
    {}                   
    { }                  
    {}                   
\theoremstyle{meiner}

\newtheorem{theorem}{Theorem}
\newtheorem{lemma}{Lemma}

\newtheorem{definition}{Definition}

\newtheorem{proposition}{Proposition}
\newtheorem{prop}{Proposition}

\newtheorem{corollary}{Corollary}

\def\Fact#1.{\par\sk{\NI\bf Fact~#1.}\ }
\def\Problem#1.{\par\sk{\NI\bf Problem~#1.}\ }
\def\Claim#1.{\medbreak\ni{\bf Claim~#1.}\ }
\def\Case#1.{\medbreak\ni{\bf Case~#1.}\ }
\def\SubCase#1.{\medbreak\ni{\bf Subcase~#1.}\ }

\def\Remark#1.{\par\sk{\NI\sffamily\bfseries Remark~#1}\ }

\def\Proof{\ni{\sl Proof.}\kern2pt\ }
\def\qed{\hfill\fbox{\hbox{}}\medskip}

\definecolor{RED}{rgb}{.84,0,0}
\definecolor{BLUE}{rgb}{0,0,.75}

\def\push(#1,#2){{\rm push}_{#1}(#2)}

\def\ovl{\overline}

\def\NN{\mathbb{N}}

\def\join{\vee}

\def\bigmeet{\bigwedge}

\def\MM{\overline{M}}

\def\X{U} 
\def\BB{\mathcal{B}}
\def\PP{\mathcal{P}}
\def\PD{\mathcal{P}_{\Delta}}
\def\BD{\mathcal{B}_{\Delta}}

\begin{document}
%

\title{\vbox{}\vskip-25mm\vbox{}
       {\bf\LARGE ULD-Lattices and $\Delta$-Bonds}}

\author{
      \large Stefan Felsner \& Kolja B.~Knauer\\[+1mm]
      \sf Institut f\"ur Mathematik,\\
      \sf Technische Universit\"at Berlin.\\
      \sf {\tt \{felsner,knauer\}@math.tu-berlin.de}
}
\date{}
\maketitle

\begin{abstract}
We provide a characterization of 
upper locally distributive lattices (ULD-lattices) in terms of
edge colorings of their cover graphs.
In many instances where a set of  combinatorial objects
carries the order structure of a lattice
this characterization yields a slick proof of distributivity or
UL-distributivity. This is exemplified by proving a 
distributive lattice structure on $\Delta$-bonds with
invariant circular flow-difference.
This instance generalizes several previously
studied lattice structures, in particular, $c$-orientations
(Propp), $\alpha$-orientations of planar graphs (Felsner, resp.~de Mendez)
and planar flows (Khuller, Naor and Klein). The characterization 
also applies to other instances, e.g.~to chip-firing games.
\end{abstract}


\section{Introduction}
\label{sec:intro}

The concept of upper locally distributive lattices (ULD) and its duals
(lower locally distributive lattices (LLD)) has appeared under several
different names, e.g. locally distributive lattices (Dilworth
\cite{Dil-60}), meet-distributive
lattices (Jamison~\cite{Jam-70,Jam-80}, Edelman~\cite{Ede-80},
Bj\"orner and Ziegler~\cite{Bjo-92b}), locally free
lattices (Nakamura~\cite{Nak-03}). Following 
Avann~\cite{Ava-68}, Monjardet~\cite{Mon-90},
Stern~\cite{Ste-99} and others we
call them ULDs. The reason for the frequent reappearance of the
concept is that there are many instances of ULDs, i.e sets of
combinatorial objects that can be \textit{naturally} ordered to form
an ULD.

ULDs have first been investigated by Dilworth~\cite{Dil-40},
many different lattice theoretical characterizations
of ULDs are known. For a survey on the work until the nineties we refer to
Monjardet~\cite{Mon-90}. We use the original definition of Dilworth: 

\begin{definition}\label{def:uld}
 Let $(P,\leq)$ be a poset. $P$ is an \term{upper locally
 distributive lattice (ULD)} if $P$ is a lattice and 
 each element has a unique minimal representation
 as meet of meet-irreducibles, i.e., there is a mapping 
 $M: P \to {\cal P}(\{ m \in P : m \text{ is meet-irreducible} \})$
 with the properties\\
 \hbox{}\quad
 $\bullet$  $ x = \bigmeet  M_x $ (representation)
 \qquad
 $\bullet$\quad $ x = \bigmeet A $ implies  $M_x \subseteq A$ (minimal).
\end{definition}

\smallskip\ni
Let $D=(V,A)$ be a directed graph, an arc coloring $c$ of $D$ is an
\term{U-coloring} if for every $u,v,w\in V$ with 
$u\neq w$ and $(v,u), (v,w) \in A$ it holds:
 \Item(U$_1$) $c(v,u)\neq c(v,w)$.
 \Item(U$_2$) There is a $z\in V$ and arcs $(u,z), (w,z)$ such that 
               $c(v,u)=c(w,z)$ and $c(v,w)=c(u,z)$. (see Figure~\ref{fig:Urule2})

   \calc_figscale{100}
    \begin{figure}[htb]
    \centerline{\input{\path/Urule2.pstex_t}}
    \caption{\label{fig:Urule2}}
    \end{figure}
    VC
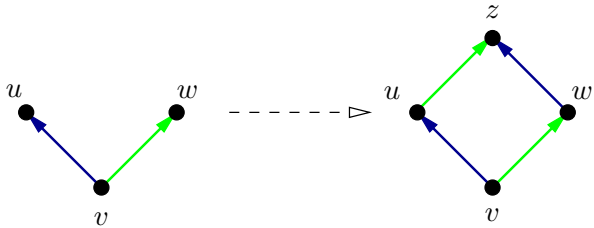
{The completion property of U-colorings.}

\begin{definition}
A finite poset $(P,\leq)$ is called \term{U-poset}
if the arcs of the cover graph $D_{P}$
of $P$ admit a U-coloring.
\end{definition}

\ni 
The characterization of ULDs in
Section~\ref{sec:characterization of ULDs}
has two parts.

\begin{theorem}\label{thm:main}
\Item(a) 
If $D$ is a finite, acyclic digraph admitting a U-coloring, then
$D$ is a cover graph, 
hence, the transitive closure of $D$ is a U-poset.\\
{\rm (b)} \
Upper locally distributive lattices are exactly the U-posets with a global minimum.
\end{theorem}

\ni
Dual in the sense of order reversal to U-coloring, U-poset and
ULD are L-coloring, L-poset and LLD.
The characterization of LLDs dual to Theorem~\ref{thm:main} 
allows easy proofs that the 
inclusion orders on the following combinatorial structures 
are lower locally distributive lattices:

\Bitem Subtrees of a tree (Boulaye~\cite{Bou-67}).
\Bitem Convex subsets of posets (Birkhoff and Bennett~\cite{Bir-85}).
\Bitem Convex subgraphs of acyclic digraphs, here a set $C$ is convex
       if $x,y\in C$ implies that all directed $(x,y)$-paths  
       are in $C$ (Pfaltz~\cite{Pfa-71}).
\medskip

\ni These combinatorial structures can also be seen as \textit{convex
sets of an abstract convex geometry}. This is no coincidence as in
fact every LLD is isomorphic to the inclusion order on the convex sets
of an abstract convex geometry and vice versa (Edelman~\cite{Ede-80}).

In Section~\ref{sec:Delta-bond} we deal with distributive lattices
arising from orientations of graphs.  To prove distributivity we use
the following well known characterization: {\em Distributive lattices
are exactly those lattices that are both ULD and LLD}
(Theorem~\ref{thm:U+L->D}).

Let $D=(V,A)$ be a connected directed graph with upper and lower
integral edge capacities $c_{u},c_{l}: A \to \mathbb{Z}$. We are
interested in maps $x: A \to \mathbb{Z}$ such that $c_{l}(a)\leq
x(a)\leq c_{u}(a)$ for all $a\in A$. The \term{circular
flow-difference} of $x$ on a cycle $C$ with a prescribed direction
is

$$
\delta(C,x):=\sum_{a\in C^{+}}x(a)-\sum_{a\in C^{-}}x(a).
$$

\ni 
Note that the circular flow-differences $\delta(C,x)$ on the
cycles of a basis of the cycle space uniquely determines the
flow-difference of $x$ on all cycles of the graph.

For a given $\Delta \in \mathbb{Z}^{\mathcal{C}}$ we consider the set
$$
\BB_{\Delta}(D,c_{l},c_{u}):=\{c_{l}\leq x\leq c_{u}\mid
\delta(C,x)=\Delta_C \text{ for all } C \in \mathcal{C}\}
$$
this is the set of \term{$\Delta$-bonds} on
$(D,c_{l},c_{u})$.  We introduce an order on $\Delta$-bonds with
prescribed circular flow-difference, i.e., on the elements of
$\BB_{\Delta}(D,c_{l},c_{u})$ such that:

\begin{theorem}\label{thm:Delta-bond-lattice}
$\BB_{\Delta}(D,c_{l},c_{u})$ carries the structure 
of a distributive lattice.
\end{theorem}

\ni
Theorem~\ref{thm:Delta-bond-lattice} is restated in 
Section~\ref{sec:Delta-bond} in a more precise version
as Theorem~\ref{thm:Delta-bond-lattice-2}. The power of this
result is exemplified in Section~\ref{sec:more-applications}
where we show that several previously studied 
distributive lattices an be recognized as special cases
of $\Delta$-bonds. This is shown in the following
cases:

\Bitem Lattice of $c$-orientations of graphs 
  (Propp~\cite{Pro-93}).
\Bitem Lattice of flow in planar graphs
  (Khuller, Naor and Klein~\cite{Khu-93}). 
\Bitem Lattice of $\alpha$-orientations of planar graphs 
  (Felsner~\cite{Fel-04}).
\medbreak

\ni
In Subsection~\ref{subsec:chip-firing} we discuss the 
\term{chip-firing game} on
directed graphs. Important properties of this game can be 
proved in the context of U-posets and ULD lattices.

\section{A Characterization of Upper Locally Distributive Lattices}
\label{sec:characterization of ULDs}

In this section we prove Theorem~\ref{thm:main}. At the end we add a
proof of the known fact that a lattice which is both ULD and LLD is
actually distributive.

The following lemma is the main tool for the proof.

\begin{lemma}\label{lem:main-lem}
Let $D=(V,A)$ be a digraph with a U-coloring $c$.
If $(x,y)$ is an arc and $p=x_0,\ldots,x_k$ a directed path from
$x=x_0$ to $z=x_k$, then there is a sequence $y_0,\ldots,y_\ell$
such that $(x_i,y_i) \in A$ and $c(x_i,y_i) = c(x,y)$ for 
$i=1,..,\ell$ and either $\ell = k$ or $\ell < k$ and
$y_\ell= x_{\ell+1}$. Figure~\ref{fig:path-chasing} illustrates the
two cases (a) and (b). Case (b) happens iff there is an edge
$(x_\ell,x_{\ell+1})$ on $p$ with $c(x,y)=c(x_\ell,x_{\ell+1})$.
\end{lemma}

   \calc_figscale{78}
    \begin{figure}[htb]
    \centerline{\input{\path/path-chasing.pstex_t}}
    \caption{\label{fig:path-chasing}}
    \end{figure}
    VC
{Illustration for Lemma~\ref{lem:main-lem}.}

\Proof Repeat applying rule U$_2$ to edges $(x_i,y_i)$ and
$(x_i,x_{i+1})$ to define a vertex $y_{i+1}$ with edges
$(y_i,y_{i+1})$ and $(x_{i+1},y_{i+1})$ such that $c(x_i,y_i) =
c(x_{i+1},y_{i+1})$.  The iteration either ends if $i=k$ (case (a)),
or if the two edges needed for the next application of the rule are
the same, i.e., $y_i=x_{i+1}$ (case (b)). In this case $c(x_i,x_{i+1})
= c(x_{i-1},y_{i-1}) = c(x,y)$, i.e., there is an edge on the path $p$
whose color equals the color of edge $(x,y)$.  Rule U$_1$ implies that
case (b) occurs whenever there is an edge on $p$ whose color equals
the color of edge $(x,y)$.  \qed

\Remark 1.
The proof does not imply that $y_i \neq x_j$ in all
cases where this is suggested by Figure~\ref{fig:path-chasing}.
An example is given Figure~\ref{fig:counterE}.
From the analysis below it follows that in all {]it bad} cases $D$ is
infinite or not acyclic.

   \calc_figscale{50}
    \begin{figure}[htb]
    \centerline{\input{\path/counterE.pstex_t}}
    \caption{\label{fig:counterE}}
    \end{figure}
    VC
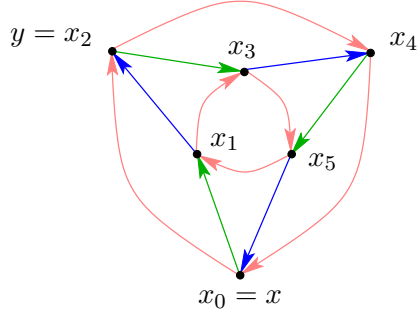
{A digraph with a U-coloring. Choosing $p = x_0,..,x_5$
and $y=x_2$ we get $y_i = x_{i+2 \pmod 6}$ for $i=0,..,6$.}

\medskip

From now on we assume that $D=(V,A)$ is a
finite, connected and acyclic digraph with a U-coloring $c$.
The assumptions imply that the transitive closure of $D$ is a finite poset
$P_D$. From the next two propositions it will follow that 
$D$ is transitively reduced, i.e., the cover graph of $P_D$.
Hence, $P_D$ is a U-poset.

\begin{proposition}\label{prop:sink}
There is a unique sink in $D$.
\end{proposition}

\Proof
Since $D$ is acyclic and finite it has a sink.
Suppose that there are two sinks $s_0$ and~$s_1$. Let $p$ be a
shortest $(s_0,s_1)$-path  in the underlying undirected graph.
In $D$ the first and the last edge of $p$ are oriented towards the sources.
Hence, there is a last vertex $x$ on $p$ which is a source on $p$,
i.e., the final part of $p$ is $y \gets x \to x_1 \to \ldots \to s_1$.
With $z = s_1$ this is the precondition of Lemma~\ref{lem:main-lem}
but both outcomes yield a contradiction: (a) is impossible because
$s_1$ is a sink, (b) is impossible because it implies that there is a
shorter  $(s_0,s_1)$-path than $p$.
\qed

We define the \term{colorset
$c(p)$ of a directed path} $p$ as the multi-set of colors 
used on edges of~$p$.

\begin{proposition}\label{prop:colorsets}
  If $D$ has a unique sink $s$ and $p$, $p'$ are directed $(x,z)$-paths, 
  then $c(p) = c(p')$. In particular this allows to define the
  \term{colorset $c(x)$ of a vertex} $x$ as the colorset of any
  $(x,s)$-path.
\end{proposition}

\Proof
First assume that $z=s$. In a top down induction we 
show that the colorsets of all $(x,s)$-paths are equal.

Assume that we know that $c(v)$ is well-defined for all vertices $v$
that are accessible from~$x$ via directed paths in $D$. Let $p$ and
$p'$ be $(x,s)$-paths with different initial edges. Let~$y$ be the
successor of $x$ on $p'$ and $a$ be the color of $(x,y)$. The edge
$(x,y)$ and the path~$p$ form the precondition for
Lemma~\ref{lem:main-lem}. Since $z=s$ is the sink the situation has to
be as in (b).  Let~$p''$ be the path $y \to y_1 \to \ldots y_\ell \to
\ldots s$.  The assumption for $y$ yields $c(p') = a \oplus c(y) = a
\oplus c(p'')$. The coloring rule U$_2$ implies that $c(x_i,x_{i+1}) =
c(y_i,y_{i+1})$ for all $i=0,..,\ell-1$ and $a = c(x,y) =
c(x_\ell,x_{\ell+1})$, therefore, $a \oplus c(p'') = c(p)$. This shows
that $c(x) = c(p') = c(p)$ is well-defined.

If $p$ is a directed $(x,z)$-path, then the concatenation $p\circ q$ 
with a $(z,s)$-path $q$ yields a $(x,s)$-path, hence,
$c(x) = c(p\circ q) = c(p) \oplus c(q) = c(p) \oplus c(z)$. This shows that
$c(p) = c(x) \ominus c(z)$ only depends on the end-vertices.
\qed

Since the colorset of a directed $(x,z)$-path  in $D$ only
depends on the end-vertices we also know that all $(x,z)$-paths
have the same length. This implies that $D$ is transitively reduced.
We have thus shown the following which is slightly stronger than statement
(a) of Theorem~\ref{thm:main}.

\begin{corollary}\label{cor:U-pos}
If $D=(V,A)$ is a
finite, connected and acyclic digraph with a U-coloring $c$,
then $D$ is a cover graph and its transitive closure is a U-poset $P_D$ with
a rank function and a $\bf 1$.
\end{corollary}

Let $P$ be a U-poset with a global minimum $\bf 0$. Define a mapping 
$\gamma : P \to \NN^k$, where $k$ is the number of colors of the U-coloring.
Assuming that the colors used by $c$ are $1,..,k$ the
$i$-th component of $\gamma(x)$ is the multiplicity of color $i$ on
any $({\bf 0},x)$-path in the cover graph $D_P$. Identifying 
vectors in $\NN^k$ with multisets of colors we have
$\gamma(x) = c({\bf 0}) - c(x)$.

\begin{lemma}\label{lem:embed}
Let $P$ be a U-poset with a global minimum $\bf 0$.
The mapping  $\gamma : P \to \NN^k$, is an order preserving embedding of 
$P$ into the dominance order on $\NN^k$.
\end{lemma}

\Proof
The implication from $y \leq_P z$ to $\gamma(y) \leq \gamma(z)$ follows
from the fact that extending a path requires more colors.
In particular the number of edges of color $i$ on a $({\bf 0},z)$-path
is at least as large as on a  $({\bf 0},y)$-path, i.e., 
$\gamma_i(y) \leq \gamma_i(z)$. In fact  $y <_P z$ 
implies $\gamma(y) \neq \gamma(z)$.

For the converse suppose $\gamma(y) \leq \gamma(z)$ but
$y \not\leq_P z$.  From the first part and $y\neq z$ we know $z \not<_P y$.
Let~$x$ be maximal with the property $x \leq_P z$ and $x\leq_P y$.
Consider the first edge $(x,y')$ on a $(x,y)$-path in $D_P$ and let
$p$ be a $(x,z)$-path. This is a situation for
Lemma~\ref{lem:main-lem}. Since the color of edge $(x,y')$ also occurs
on $p$ we are in case~(b). This case, however, is
impossible because $y' \leq_P z$,  $y'\leq_P y$ and $y' \geq x$
contradicts the choice of $x$.
\qed

\begin{lemma}\label{lem:join}
Let $P$ and $\gamma$ be as above.
For all $z,y \in P$ there is a $w\in P$ with $\gamma(w)
= \gamma(y) \join \gamma(z)$ where $\join$ is the componentwise
maximum.
\end{lemma}

\Proof
For any fixed $y$ we proceed with top down induction. Given $z$
consider a  maximal $x$ with the property $x \leq_P z$ and $x\leq_P y$.
Let $(x,y')$ be the first edge on a $(x,y)$-path in $D_P$ and 
let $p$ be a $(x,z)$-path. Case (b) of Lemma~\ref{lem:main-lem} is
impossible because $y'$ would have prevented us from choosing $x$.
Hence, we are in case (a) and there is a $z'$ covering $z$ such that 
the edges $(z,z')$ and $(x,y')$ have the same color $i$, moreover,
the path from $x$ to $z$ has no edge of color $i$. Induction implies that
there is a  $w'$ such that  $\gamma(w') = \gamma(z') \join \gamma(y)$.
Since $\gamma(z') = \gamma(z) + {\bf e}_i$ and
for the $i$-th component  $\gamma_i(z) = \gamma_i(x) < \gamma_i(y)$
holds we can conclude $\gamma(z') \join \gamma(y) = \gamma(z) \join \gamma(y)$,
i.e., $w'$ may also serve as $w$.
\qed

\begin{proposition}\label{prop:lattice}
If a U-poset has a global minimum $\bf 0$, then it is a lattice.
\end{proposition}

\Proof The mapping $\gamma$ is an order embedding of $P$ into $\NN^k$
(Lemma~\ref{lem:embed}) and the image $\gamma(P)$ is join-closed
(Lemma~\ref{lem:join}).  Together this implies that $P$ has unique
least upper covers (joins exist). Since $P$ has a $\bf 0$ there is a
lower cover for every pair of elements. Uniqueness for greatest lower
covers (meets) follows from the existence of unique joins.  \qed

The next goal is to show that every element of $P$ has a unique minimal
representation as a meet of meet-irreducibles. Let $C(x)$ be the set of
colors of the edges emanating from $x$ in the cover graph. With the next lemma
we associate a meet irreducible element with every color $i\in C(x)$.

\begin{lemma}\label{lem:mi-cands}
  Let $P$ be a U-poset with a $\bf 0$ and let $x\neq {\bf 1}$ be an
  element of $P$.  For every $i\in C(x)$ there is a unique maximal
  element $y_i$ such that
\quad$\bullet$ 
  $y_i \geq x$ and 
\quad$\bullet$ 
  $\gamma_i(y_i) = \gamma_i(x)$.\\
  The element $y_i$ is meet-irreducible and $\gamma_j(y_i) >
  \gamma_j(x)$ for all $j\in C(x)\setminus\{i\}$.
\end{lemma}

\Proof Let $i\in C(x)$ and consider the set $S_i(x)$ of all $y \geq x$ 
with $\gamma_i(y) = \gamma_i(x)$.  The set $S_i(x)$ contains $x$, 
hence, it is non-empty and by Lemma~\ref{lem:join} it contains a unique maximal
element~$y_i$. The element $y_i$ is meet-irreducible, otherwise
we could find a successor of $y_i$ in $S_i(x)$. For every $j\in C(x)$
there is an element $x_j$ with $\gamma(x_j) = \gamma(x) + {\bf e}_j$.
For $j\neq i$ the element $x_j$ is in $S_i(x)$, hence, 
$\gamma(x_j) \leq \gamma(y_i)$ and $\gamma_j(x) < \gamma_j(x_j) \leq \gamma_j(y_i)$.
\qed

\begin{proposition}\label{prop:uld}
A U-poset with a global minimum  $\bf 0$ is an upper locally distributive lattice.
\end{proposition}

\Proof
We claim that $M_x = \{ y_i : i \in C(x) \}$ is the unique
minimal set of meet-irreducibles with $x = \bigmeet M_x$.  

Let $z$ be any lower bound for $M_x$, i.e., an
element with $z \leq y_i$ for all $i\in C(x)$.
Since $\gamma$ is order preserving and $\NN^k$ is closed under
taking meets we have $\gamma(z) \leq \bigmeet \{\gamma(y_i) : i\in C(x) \}$.
From $\gamma_i(y_i) = \gamma_i(x)$ it follows that
$\bigmeet \{\gamma(y_i) : i\in C(x) \} = \gamma(x)$. 
Since $\gamma$ is order preserving this implies $z \leq x$,
i.e., $x$ is the unique maximal lower bound for $M_x$ 
and the notation $x = \bigmeet M_x$ is justified.

It remains to show that the representation $x = \bigmeet M_x$ is the
unique minimal representation of $x$ as meet of meet-irreducibles.
Let $i\in C(x)$ and consider a set $M$ of meet-irreducibles with $y_i
\not\in M$. It is enough to show that $x \neq \bigmeet M$. If $M$
contains a $y$ with $x \not< y$, then $x \neq \bigmeet M$ is obvious.
Consider the set $S_i(x)$ from the proof of Lemma~\ref{lem:mi-cands},
every element $y\neq y_i$ in this set is contained in a $(x,y_i)$-path
$p$ that contains no $i$-colored edge. Lemma~\ref{lem:main-lem}
implies that there is an $i$-colored edge leaving $y$ together with
the edge leaving $y$ on $p$ this implies that $y$ is not
meet-irreducible. Hence $M\cap S_i(x) = \emptyset$. All
$y > x$ with $y\not\in S_i(x)$ satisfy $y \geq x_i$ with our
assumptions on $M$ this implies that $x_i$ 
with $\gamma(x_i) = \gamma(x) + {\bf e}_i$ is a lower bound on $M$,
i.e., $x \neq \bigmeet M$.
\qed

From what we have shown so far we obtain the following criterion
for ULD lattices.

\begin{theorem}\label{thm:criterion}
If $D$ is a finite, acyclic digraph with a unique source ${\bf 0}$
and there is a U-coloring of the arcs of $D$, then the transitive 
closure of $D$ is an upper locally distributive lattice.
\end{theorem}

To complete the proof of Theorem~\ref{thm:main} it remains to 
show that every ULD has a representation as U-poset, i.e., 
we have to present a U-coloring of its cover graph.

\begin{theorem}\label{thm:uld-is-u}
The cover graph of every finite ULD lattice admits a U-coloring.
\end{theorem}

\ni Consider the mapping $\MM$ that takes an element $x$ of $P$ to the
set $\MM_x$ of all meet-irreducible elements that are at least as
large as $x$.  The definition of meet irreducible implies that $x =
\bigmeet\MM_x$ for all $x$, i.e., the set $\MM_x$ uniquely
characterizes $x$. Moreover, $x \leq y$ iff $\MM_x \supseteq \MM_y$.

On the basis of the mappings $M$ and $\MM$ we will define a
U-coloring of the cover relations of $P$.  As colors we use the
meet-irreducible elements of $P$.

\begin{lemma}\label{lem:cover}
Let $P$ be a ULD lattice.
A comparability $x < y$ is a cover  iff \ $|\MM_x\setminus \MM_y| = 1$.
\end{lemma}

\Proof 
An element $z$ with $x < z < y$ satisfies $\MM_y \subsetneq
\MM_z \subsetneq \MM_x$ which implies $|\MM_x\setminus \MM_y| \geq 2$.

Let $x < y$ and suppose that $|\MM_x\setminus \MM_y| \geq 2$.  Since
$\bigmeet M_x < \bigmeet \MM_y$ there has to be some $m\in
M_x\setminus \MM_y$.  Let $z = \bigmeet(\MM_x - m)$, we claim that
$x < y < z$. From $\MM_x - m \subset \MM_x$ it
follows that $z \geq x$, and since $m\in M_x$ and $m\not\in M_z$ 
Definition~\ref{def:uld} implies that $z < x$. 
Since $(\MM_x - m) \supseteq \MM_y$ we have
$z \leq y$. Let $m'$ be an element with $m \neq m' \in \MM_x\setminus
\MM_y$, it follows that $m' \in \MM_z$ and $m' \not\in \MM_y$.
Therefore $z\neq y$ and we have shown that $x < z < y$, i.e., the pair
$x,y$ is not in a cover relation.  \qed

To a cover relation $x \prec y$ we assign the unique meet-irreducible in
$\MM_x\setminus \MM_y$ as its color. Note that this meet-irreducible
is a member of $M_x$.
To verify that this is a U-coloring we have to
check the two properties U$_1$ and U$_2$. 

\Claim 1. The coloring satisfies U$_1$.

\Proof
Let $x \prec y_1$ and $x\prec y_2$ be two cover relations. 
Since $x$ is the meet of $y_1$ and $y_2$ we have
the representation $x = \bigmeet(\MM_{y_1} \cup \MM_{y_2})$ of $x$ as
meet of irreducibles, hence,
$M_x \subseteq \MM_{y_1} \cup \MM_{y_2}$.
If both covers had the same color $m$, then $m\in M_x$ 
but $m\not\in \MM_{y_1} \cup \MM_{y_2}$, a contradiction.
\qed

\Claim 2. The coloring satisfies U$_2$.

\Proof
Let $x \prec y_1$ and $x\prec y_2$ be two cover relations
such that $x \prec y_1$ has color $m_1$
and $x\prec y_2$ has color $m_2$, i.e., $\MM_{y_i} = \MM_x - m_i$. 
Consider $z = \bigmeet( \MM_x - m_1 - m_2)$. Since $z$ is representable as
meet of elements from $\MM_{y_i}$ we know $z \geq y_i$ for $i=1,2$.
Since $y_1$ and $y_2$ both cover $x$ it follows that $z \neq y_1,y_2$,
hence $z > y_i$ for $i=1,2$. From $\MM_x - m_1 - m_2 \subseteq \MM_z 
\subset \MM_{y_i} = \MM_x - m_i$ it follows that 
$|\MM_{y_i}\setminus \MM_z| = 1$. Lemma~\ref{lem:cover} implies that
$z$ covers each  $y_i$ and the labels of these covers are as required.
\qed

In many applications of the characterization of ULDs 
the lattice in question is actually distributive.
Such a situation is the topic of the next section.
To make the paper self contained we prove the following
folklore result.

\begin{theorem}\label{thm:U+L->D}
If a finite, acyclic and connected digraph $D$ 
admits a U- and a L-coloring
then $D$ is the cover graph of a distributive lattice $P_D$.
Moreover, the colorings yield an explicit cover-preserving
embedding $P_D \hookrightarrow \NN^k$, where $k$ is the number of colors.
\end{theorem}

\Proof
Corollary~\ref{cor:U-pos} and its
dual imply that $D$ is the cover graph of a
poset $P_D$ with $\bf 0$ and  $\bf 1$. Hence,
with Proposition~\ref{prop:lattice} $P_D$ is
a lattice.

Let $c_U$ and $c_L$ be a U- and a L-coloring of $D$.
Consider the coloring $c = c_U \times c_L$. 
The claim is that $c$ is both a U- and a L-coloring of $D$.
The rule U$_1$ and its dual L$_1$ are immediately inherited
from the corresponding rules for  $c_U$ and $c_L$.
Whenever there is a diamond $x \prec y_1$, $x\prec y_2$,
$y_1 \prec z$, $y_2 \prec z$ the colors of a pair of
$(x,y_i)$ and $(y_j,z)$ in $c_U$ and $c_L$ coincides.
This implies rules U$_2$ and L$_2$ for $c$.

Consider the order embedding $\gamma : P_D \to \NN^k$ 
that is based on the coloring $c$. By
Lemma~\ref{lem:join}~$\gamma$ is compatible with joins, the dual
implies that $\gamma$ is compatible with meets. 
Therefore~$\gamma$ is a lattice embedding, i.e., $P$ a
sublattice of the distributive lattice $\NN^k$,
hence, $P$ is itself distributive.
\qed

\begin{wrapfigure}[7]{R}{0pt}
   \calc_figscale{50}
   \input{\path/semi-uld.pstex_t}
    
\caption{\label{fig:semi-uld}}
\end{wrapfigure}

\NI\sffamily\bfseries Remark~2 \rm
Let $D$ be a digraph with a U-coloring. 
We need acyclicity, connectedness, finiteness
and the unique source to conclude that $D$ 
corresponds to a finite ULD lattice. We feel that among these
conditions the {\it unique source} has a somewhat
artificial flavour. 
Abstaining on this condition
it can be shown (along the lines of our proof) 
that the corresponding poset $P$
is a join-semilattice with the property
that for all $x\in P$ there is a unique minimal
set $M_x$ of meet-irreducibles such that
$x$ is a maximal lower bound for $M_x$.
Figure~\ref{fig:semi-uld} shows a small example,
in this case $M_s=M_t=\{u,v\}$.

\section{The Lattice of $\Delta$-Bonds}
\label{sec:Delta-bond}

Recall the setting from the introduction: The data are a
directed multi-graph $D=(V,A)$ with upper and lower integral edge
capacities $c_{u},c_{l}: A \to \mathbb{Z}$ and a number $\Delta_C$ for
each cycle $C$ of $D$.
We are interested in the \term{$\Delta$-bond}\ncterm{s} for this data, i.e,
in maps $x: A \to \mathbb{Z}$ such that 

\Item(D$_1$) \quad $c_{l}(a)\leq x(a)\leq c_{u}(a)$ for all $a\in A$.
             \hfill (capacity constraints)

\Item(D$_2$) \quad $\displaystyle \Delta_C = \sum_{a\in C^{+}}x(a)-
             \sum_{a\in C^{-}}x(a) = \delta(C,x)$ for all $C$.
             \hfill (circular flow conditions)

\ni
Throughout the discussion we shall assume that the data
$(D,c_{l},c_{u},\Delta)$ are such that the set of
corresponding $\Delta$-bonds is non-empty. Moreover, we
want to simplify matters by concentrating on connected
graphs and getting rid of \term{rigid} edges, these are
edges $a\in A$ with $x(a) = y(a)$ for all pairs $x$, $y$
of $\Delta$-bonds.

Let $a$ be a rigid edge of $D$ and let $D/_a$ be obtained from $D$ by contracting
edge $a$. Since we allow multiple edges and loops the cycles in $D/_a$ and in $D$ are
in bijection. Let $C/_a$ be the cycle in $D/_a$ corresponding to $C$ in $D$.
Define $\Delta'_{C/_a} = \Delta_C$ if $a\not\in C$ and 
$\Delta'_{C/_a} = \Delta_C - x(a)$ if $a\in C^+$  and 
$\Delta'_{C/_a} = \Delta_C + x(a)$ if $a\in C^-$.
These settings yield the bijection that proves
\begin{lemma}\label{lem:rigid}
$\BB_{\Delta}(D,c_{l},c_{u}) \cong \BB_{\Delta'}(D/_a,c_{l},c_{u})$.
\end{lemma}

Given data $(D^1,c^1_{l},c^1_{u},\Delta^1)$ and $(D^2,c^2_{l},c^2_{u},\Delta^2)$
there is an obvious extension to a union structure $(D,c_{l},c_{u},\Delta)$ where
$D$ is the union of graphs and the $c_{l},c_{u},\Delta$ are concatenations
of vectors. Since $\Delta$-bonds factor into a $\Delta_1$- and a $\Delta_2$-bond
we have:
\begin{lemma}\label{lem:prod}
$\BB_{\Delta}(D,c_{l},c_{u}) \cong 
\BB_{\Delta^1}(D^1,c^1_{l},c^1_{u}) \times \BB_{\Delta^2}(D^2,c^2_{l},c^2_{u})$.
\end{lemma}

The data $(D,c_{l},c_{u},\Delta)$ are \term{reduced} if $D$ is connected and 
there is no rigid edge. Henceforth we will assume that any given set of data is
reduced.

With a partition $(\X,\ovl{\X})$ of the vertices $V$ of $D$ we consider
the cut $S=S{[\X]} \subset A$. The \ncterm{forward edges} $S^+$ of $S$
are those $a\in A$ directed from $\X$ to $\ovl{\X}$, \ncterm{backward
  edges} $S^-$ of $S$ are directed from $\ovl{\X}$ to $\X$. 

For $x: A \to \mathbb{Z}$ and a subset $\X \subset V$ we define
$y=\push(U,x)$ such that $y(a) = x(a) +1$ for all $a\in S^+[\X]$,
$y(a) = x(a) -1$ for all $a\in S^-[\X]$ and $y(a) = x(a)$ for all
$a\not\in S[\X]$. We say that $y=\push(U,x)$ is obtained by
\term{pushing $\X$} in $x$.

Fix an arbitrary vertex $v_0$ in $D$ as the \term{forbidden vertex}.
For $x,y \in\BD(D,c_{l},c_{u})$ define $x \leq y$ if $y$ can be reached
from $x$ via a sequence of pushes at sets $\X_i$, such that
$v_0\not\in \X_i$ for all~$i$ and all intermediate states in the
sequence are in $\BD(D,c_{l},c_{u})$ .

Below (Corollary~\ref{cor:order}) we show that the relation ``$\leq$'' makes
$\BD(D,c_{l},c_{u})$ into a partial order $\PD(D,c_{l},c_{u})$.

The main result of this section is Theorem~\ref{thm:Delta-bond-lattice}
which can now be stated more precisely:

\begin{theorem}\label{thm:Delta-bond-lattice-2}
The order $\mathcal{P}_{\Delta}(D,c_{l},c_{u})$ on $\Delta$-bonds 
is a distributive lattice.
\end{theorem}

It would have been be more precise to write
$\mathcal{P}_{\Delta}(D,c_{l},c_{u},v_0)$ in the theorem
because the actual lattice depends on the choice of the forbidden vertex.
Different forbidden vertices yield different lattices on the same ground set.
The result of Lemma~\ref{lem:prod} carries over to the lattices, i.e.,
if $D^1$ and $D^2$ are connected, then
$\PD(D,c_{l},c_{u},\{v_1,v_2\}) \cong
\PP_{\Delta^1}(D^1,c^1_{l},c^1_{u},v_1) \times
\PP_{\Delta^2}(D^2,c^2_{l},c^2_{u},v_2)$.

The next lemma describes the condition for a legal push, i.e., a push that
transforms a $\Delta$-bond into a $\Delta$-bond.

\begin{lemma}\label{lem:legal}
If $x$ is in $\BB_{\Delta}(D,c_{l},c_{u})$ 
and $S=S{[\X]}$ is a cut such that $x(a) < c_u(a)$ for all $a\in S^+$
and $x(a) > c_l(a)$ for all $a\in S^-$, then $y=\push(U,x)$
is also in $\BB_{\Delta}(D,c_{l},c_{u})$.
\end{lemma}

\proof The assumption on $c_u(a)$ and $c_l(a)$ implies that $y=\push(U,x)$
respects the capacity constraints. From the orthogonality of
the cycle space and the bond space of $D$ it follows that
$\delta(C,x) = \delta(C,y)$ for all cycles $C$, i.e., $y$ 
satisfies the circular flow conditions.
\qed

\begin{lemma}\label{lem:v-push}
If $x \leq y$, then $y$ can  be
obtained from $x$ by a sequence of \term{vertex push}\ncterm{s},
i.e., at cuts $S[v]$ with $v\neq v_0$. 
\end{lemma}

\proof It is enough to show that a single push of a set $\X$ can be
replaced by a sequence of vertex pushes.  Let $y=\push(U,x)$ we show
that there is a vertex $w \in \X$ such that $x'=\push(w,x)$
is a $\Delta$-bond and $y=\push(\X\setminus w,x')$.
This implies the result via induction on the size of $U$.

Choose $w \in \X$ arbitrarily. If pushing $w$ is legal,
i.e., $x(a) < c_u(a)$ for all $a\in S[w]^+$
and $x(a) > c_l(a)$ for all $a\in S[w]^-$, then let $x'=\push(w,x)$.
Note that pushing $\X\setminus w$ in  $x'$ is legal:
Indeed for $a\in S[\X\setminus w] \cap S[\X]$ we have $x'(a) = x(a)$
and all $a\in S[\X\setminus w] \setminus S[\X]$ are incident to $w$,
if such an $a$ is in $S[\X\setminus w]^+$, then it is in $S[w]^-$ and
$x'(a) = x(a)-1 \leq c_u(a) -1$, hence, $x'(a) < c_u(a)$. 
The case for $a\in S[\X\setminus w]^-$ is symmetric. 

It remains to show that there is a $w\in \X$ such that pushing $w$ 
is legal. Choose $w_1\in \X$ arbitrarily if
$S[w_1]$ is not legal, then there is an incident edge $a_1\in
S[w_1]^+$ with $x(a_1) = c_u(a_1)$ or $a_1\in S[w_1]^-$ with $x(a_1) =
c_l(a_1)$. Let $w_2$ be the second vertex of $a_1$ and note that
$w_2\in \X$. If pushing $w_2$ is not legal, then there is an incident
edge $a_2$ obstructing the push and so on. This yields a
sequence $w_1,a_1,w_2,a_2,\ldots$, either the sequence ends in a
vertex $w_j$ which is legal for pushing or it closes into a cycle.
Assume that there is a cycle $C$ such that $x(a) = c_u(a)$ for all $a
\in C^+$ and $x(a) = c_l(a)$ for all $a \in C^-$.  The condition
implies $\Delta_C = \delta(C,x) = \sum_{a\in C^{+}}c_u(a)-\sum_{a\in
  C^{-}}c_l(a)$.  It follows that every $\Delta$-bond $y$ has $y(a) =
x(a)$ for all $a\in C$, i.e., the edges in $C$ are rigid. 
However, the data $(D,c_u,c_l,\Delta)$ are assumed to be reduced, i.e.,
there are no rigid edges.  Hence, there must be a $w\in \X$ 
that is legal for pushing.  \qed

In the case of a general set of data $(D,c_u,c_l,\Delta)$ with rigid
edges we would have to allow the pushing of $\X$ iff $\X$ is the
vertex set of a connected component of rigid edges.

\begin{corollary}\label{cor:order}
The relation $\leq$ on $\BB_{\Delta}(D,c_{l},c_{u})$ is acyclic,
i.e., it is an order relation. 
\end{corollary}
\proof 
Otherwise we could linearly combine
vertex cuts $S[v]$, with $v\neq v_0$, to zero. But these vertex cuts
are a basis of the bond space.\qed
 
Pushes at vertices correspond to the cover relations, i.e., edges of
the cover graph, of $\mathcal{P}_{\Delta}(D,c_{l},c_{u})$
(Lemma~\ref{lem:v-push}).  A coloring of the edges of the cover
graph of $\PD = \mathcal{P}_{\Delta}(D,c_{l},c_{u})$ with colors
from $V\setminus\{v_0\}$ is naturally given by coloring a cover with
the vertex of the corresponding vertex cut.

\begin{lemma}
The coloring of the edges of the cover graph of
$\mathcal{P}_{\Delta}$ with colors in $V\setminus\{v_0\}$
is a U-coloring.
\end{lemma}

\Proof 
Let $x\in \BB_{\Delta}(D,c_{l},c_{u})$ and suppose that pushing $v$ in $x$
is legal, i.e., there is a covering colored $v$
leaving $x$. The other element $x'$ of the covering pair is completely
determined by $x$ and $v$. This shows property U$_1$.

For U$_2$ assume that $u$ and $v$ can both be pushed in $x$. We
have to show that they can be pushed in either order. This
clearly holds if the vertex cuts of $u$ and $v$ are disjoint.
Suppose $(u,v)\in A$ and note that $c_{l}(a)<x(a)<c_{u}$. Now after
pushing $u$ we still have $c_{l}(a)<x(a)$ thus we can still
push at $v$.  Conversely pushing $v$ preserves $x(a)<c_{u}$,
i.e., the push of $u$ remains legal.  \qed

A completely symmetric argument shows that the coloring is also
a U-coloring for the reversed order, i.e., a L-coloring.
Theorem~\ref{thm:U+L->D} implies that every connected component
of $\PD$ is a distributive lattice. To complete the proof of
Theorem~\ref{thm:Delta-bond-lattice} it only remains to show that 
$\PD$ is connected. This is shown in the last lemma of this section.

\begin{lemma}
The order  $\PD$ is connected.
\end{lemma}

\proof
With $\Delta$-bonds $x$ and $y$ consider $z=x-y$. 
Note that $z$ is a $\bf 0$-bond because 
$\delta(C,z) = \delta(C,x) - \delta(C,y) = 0$ for all cycles $C$.
Since $\bf 0$-bonds are just bonds there is a unique expression
of $z$ as a linear combination of vertex cuts $S[v]$ with $v\neq v_0$,
we write this as $z = \sum_v \lambda_v S[v]$.
based on the coefficients $\lambda_v$ define 
$T =\{ v\in V : \lambda_v > 0\}$. Since $x$ and $y$ 
are different and could as well be 
exchanged we may assume that $T\neq\emptyset$.

We claim that pushing $T$ in $y$ is legal, i.e., $y' = \push(T,y)$ is
a $\Delta$-bond. First note that $v_0 \not\in T$. Now let $a = (v,w)$
be an arc with $v\in T$ and $w \in V\setminus T$, i.e., $a \in
S[T]^+$, from $(x-y)(a) = \lambda_v - \lambda_u >0$ we obtain $x(a) >
y(a)$, hence, $y(a) < c_u(a)$. For $a \in S[T]^-$ we obtain $x(a) <
y(a)$, hence, $y(a) > c_l(a)$.  Lemma~\ref{lem:legal} implies that
$y'$ is a $\Delta$-bond.

Assuming that $\PD$ has several components we may choose $x$ and $y$
from different components such that $\sum_v \lambda_v$ is minimal,
where $x-y = \sum_v \lambda_v S[v]$. Since $y'$ is obtained from $y$
by pushing the set $T$ it is in the same component as $y$. Note that $x -
y' = \sum_{v\in V\setminus T} \lambda_v S[v] + \sum_{v\in T}
(\lambda_v-1) S[v]$ since $\sum_{v} \lambda_v - |T| <
\sum_v \lambda_v$ this contradicts the choice of $x$ and $y$.
\qed

\section{More Applications of the ULD Characterization}
\label{sec:more-applications}

In the first three parts of this section we deal with special cases of
Theorem~\ref{thm:Delta-bond-lattice-2}. As a result we reprove known
instances of distributive lattices from graphs.
Subsection~\ref{subsec:potentials} connects from $\Delta$-bonds, i.e.,
special edge weightings, to potentials, i.e., vertex weightings.
In Subsection~\ref{subsec:chip-firing} we discuss the
chip-firing game on directed graphs. A central and previously
know result is that the states of this game carry the structure of an
ULD lattice. We obtain this as a direct application of our
characterization.

\subsection{The lattice of $c$-orientations -- Propp~\cite{Pro-93}}

Given an orientation $O$ of a graph $G=(V,E)$. We regard a cycle of
an undirected graph as an edge-sequence, rather than an edge-set, i.e.
a cycle comes with the direction of its traversal.  If $C$ is a cycle
in $G$ then denote by $c_O(C):=|C_O^+| - |C_O^-|$ the \term{circular
  flow-difference} of $O$ around $C$, where $C_O^+$ is the set of
forward arcs of $C$ in $O$ and $C_O^-$ is the set of backward arcs.

Given a vector $c$, which assigns to every cycle $C$ of $G$ an integer
$c(C)$, we call an orientation $O$ of $G$ with $c(C) = c_O(C)$ a
\term{$c$-orientation}.

The main result in Propp's article~\cite{Pro-93} is:
\begin{theorem}\label{thm:pro}
Let $G=(V,A)$ be a graph and $c\in\mathbb{Z}^\mathcal{C}$.
The set of $c$-orientations of $G$ carries the structure
of a distributive lattice.
\end{theorem}

\begin{proof}
  Let $D=(V,A)$ be any orientation of $G$. Define
  $\Delta:=\frac{1}{2}(c_D-c)$. We interpret
  $x\in\mathcal{B}_{\Delta}(D,\mathbf{0},\mathbf{1})$ as the
  orientation $O(x)$ of $G$ which arises from $D$ by changing the
  orientation of $a\in A$ if $x(a)=1$. For an arc set $A'\subseteq A$
  denote by $x(A'):=\sum_{a\in A'}x(a)$. We calculate

\begin{eqnarray*}
	c_{O(x)}(C)     &=&|C_{O(x)}^{+}|-|C_{O(x)}^{-}|\\
		        &=&|C_D^+|-x(C_D^+)+x(C_D^-)-(|C_D^-|-x(C_D^-)+x(C_D^+))\\
			&=&|C_D^+|-|C_D^-|-2(x(C_D^+)-x(C_D^-))\\
			&=&c_D(C)-2\delta(x,C) = c_D(C)-2\Delta_C\\
			&=&c
\end{eqnarray*}

This shows that $c$-orientations of $G$ correspond bijectively to
$\Delta$-bonds in $\mathcal{B}_{\Delta}(D,\mathbf{0},\mathbf{1})$. By
Theorem~\ref{thm:Delta-bond-lattice-2} we obtain a distributive
lattice structure on the set of $c$-orientations of $G$.
\end{proof}

From duality of planar graphs and the above theorem Propp
derives the following two corollaries:

{\Bitem The set of $d$-factors of a plane bipartite graph can be
        enhanced with a distributive lattice structure.
\Bitem The set of spanning trees of a  plane graph can be
        enhanced with a distributive lattice structure.

}

\subsection{The lattice of flow in planar graphs  -- Khuller, Naor
and Klein~\cite{Khu-93}}

Consider a planar digraph $D=(V,A)$, with each arc $a$ having an
integer \nct{lower} and \nct{upper} bound on its capacity, denoted
$c_{\ell}(a)$ and $c_{u}(a)$.  For a function $f:A\to \mathbb{Z}$ call
$\omega(v,f):=\sum_{a\in {\rm in}(v)} f(a) - \sum_{a\in {\rm out}(v)}
f(a)$ the \term{excess} at $v$. Given a vector
$\Omega\in\mathbb{N}^{V}$ call $f$ a \term{$\Omega$-flow} if
$c_{\ell}(a) \leq f(a)\leq c_{u}(a)$ for all $a$ and
$\Omega_v=\omega(v,f)$ for all $v\in V$. Denote by
$\mathcal{F}_{\Omega}(D,c_{\ell},c_u)$ the set of $\Omega$-flows.

\begin{theorem}\label{thm:genkhu}
  If $D$ is a planar digraph then
  $\mathcal{F}_{\Omega}(D,c_{\ell},c_u)$ carries the structure of a
  distributive lattice.
\end{theorem}

\begin{proof} Given a crossing-free plane embedding of $D$ we look at
  the planar dual digraph $D^{*}$. It is an orientation of the planar
  dual $G^{*}$ of the underlying undirected graph $G$ of $D$. Let $v$
  be a vertex of $G^{*}$ corresponding to a facial cycle $C$ of the
  embedding of $D$. Orient an edge incident to $v$ as outgoing arc of
  $v$ if the dual arc is forward when traversing $C$ in clockwise
  direction.  Given values on the arcs of $D$ we simply transfer them
  to the corresponding arcs of $D^{*}$.

  Since the excess at a vertex of $D$ dualizes to the circular flow
  difference around the corresponding facial cycle of
  $D^{*}$\kern-2pt, we have a correspondence between
  $\mathcal{F}_{\Omega}(D,c_{\ell},c_u)$ and
  $\mathcal{B}_{\Omega}(D^{*},c_{\ell},c_u)$\kern-1pt. This yields the
  distributive lattice structure on $\Omega$-flows of planar graphs.
\end{proof}

Analogous to the case of $\Delta$-bonds we can assume the data
$(D,c_{\ell},c_u,\Omega)$ to be reduced. Now the dual operation to
vertex pushes is to augment the flow around facial cycles. A natural
candidate for the \term{forbidden facial cycle} is the unbounded face
of the planar embedding. By flow-augmentation at the remaining facial
cycles we can construct the cover graph of a distributive lattice on
$\mathcal{F}_{\Omega}(D,c_{\ell},c_u)$.

Khuller, Naor and Klein~\cite{Khu-93}, actually, only consider the
special case of Theorem~\ref{thm:genkhu} where $\Omega = \mathbf{0}$,
these flows without excess are called \term{circulation}\nct{s}.  We
restate their result:

\begin{theorem}
  Let $D$ be a planar digraph with upper and lower arc capacities
  $c_l$ and $c_u$. The set of circulations of $D$ within $c_l$ and $c_u$
  carries the structure of a distributive lattice.
\end{theorem}
	
\subsection{The lattice of $\alpha$-orientations in planar graphs  -- 
Felsner~\cite{Fel-04}}

Consider a plane graph $G=(V,E)$. Given a mapping $\alpha: V \to
\mathbb{N}$ an orientation $X$ of the edges of $G$ is called an
\term{\hbox{$\alpha$-orientation}} if $\alpha$ records the out-degrees
of all vertices, i.e.,${\sf outdeg}_X(v) = \alpha(v)$ for all $v\in V$.

The main result in~\cite{Fel-04} is:

\begin{theorem}\label{thm:alpha}
  Given a planar graph and a mapping $\alpha: V \to \mathbb{N}$ the
  set of $\alpha$-orientations of $G$ carries the structure of a
  distributive lattice.
\end{theorem}
\begin{proof}
  Analogously to the proof of Theorem~\ref{thm:pro}, where
  c-orientations were interpreted as elements of
  $\mathcal{B}_{\Delta}(C,\mathbf{0},\mathbf{1})$ we can view
  $\alpha$-orientations as elements of
  $\mathcal{F}_{\Omega}(D,\mathbf{0},\mathbf{1})$ for some orientation
  $D$ of $G$ and $\Omega$ depending on $\alpha$. Application of
  Theorem~\ref{thm:genkhu} yields the distributive lattice structure
  on the $\alpha$-orientations.
\end{proof}

In view of the present paper $\alpha$-orientations appear as a special
case of the preceeding constructions. Nevertheless they are quite
general objects. Special instances of $\alpha$-orientations on plane
graphs yield lattice structures on

{\Bitem
Eulerian orientations of a plane graph.
\Bitem
Spanning trees and $d$-factors of a plane graph.
\Bitem
Schnyder woods of a 3-connected plane graph.
}
\medskip

\ni It is natural to ask for the structure of $\alpha$-orientations of
graphs embedded on some surface. Propp~\cite{Pro-93} comments that to
move between the $d$-factors in toroidal graphs it is necessary to
operate on non-contractible cycles. In his diploma-thesis
Knauer~\cite{Kna-07} investigates $\alpha$-orientations on
surfaces and generalizations of $\alpha$-orientations.

\subsection{Potentials}
\label{subsec:potentials}

For a last and fundamental example set $\Delta=\mathbf{0}$. Fix
$v_0\in V$ and denote by $x(A'):=\sum_{a\in A'}x(a)$ where
$A'\subseteq A$. Now $x\in \mathcal{B}_{\mathbf{0}}(D,c_{\ell},c_u)$
means that the potential mapping $\pi_{x}(v):=x(p)$ for a
$(v_{0},v)$-path~$p$ is well defined, i.e., independent of the choice
of~$p$. In particular for every $v_{0}\in V$ the set
$\mathcal{B}_{\mathbf{0}}(D,c_{\ell},c_u)$ is in bijection with the
set of \term{feasible vertex-potentials}:

$$
\Pi_{v_{0}}(D,c_{\ell},c_u):=
  \{\pi\in\mathbb{N}^{V}\mid\pi(v_{0})=0 \text{ and }
  c_{\ell}(a)\leq\pi(w)-\pi(w)\leq c_{u}(a) \text{ for all } a\in A\}.
$$

\ni
Theorem~\ref{thm:Delta-bond-lattice-2} yields that
$\Pi_{v_{0}}(D,c_{\ell},c_u)$ carries the structure of a distributive
lattice. Indeed this can be obtained easier since
$\Pi_{v_{0}}(D,c_{\ell},c_u)$ is a suborder of the
distributive lattice $\mathbb{N}^{V}$ where the order is given by
dominance and the lattice operations are max and min.
The simple observation that for $\pi, \pi'\in
\Pi_{v_{0}}(D,c_{\ell},c_u)$ also their componentwise minimum and
maximum are in $\Pi_{v_{0}}(D,c_{\ell},c_u)$ yields the distributive
lattice structure on feasible vertex potentials.

The proof given in~\cite{Fel-04,Khu-93,Pro-93} for the respective
lattice structures are all based on the construction of potentials
corresponding to the objects they investigate.  In a forthcoming
paper~\cite{FeKn-08} we exploit the potential approach. In this paper
will deal with a class of {\it distributive polytopes},
\term{{\rm D}-polytopes} for short and we generalize the notion of
$\Delta$-bonds to \textit{generalized bonds}. These are in a certain
sense the most general structure on graphs that form a distributive
lattice.

\subsection{Chip-Firing Games}
\label{subsec:chip-firing}

A well-studied class of ULDs come from \textit{chip-firing games}. As 
in~\cite{Bjo-92a} we consider chip-firing games on directed graphs. For
the undirected case see~\cite{Bjo-91}. In general the ULDs coming
from chip-firing games are not distributive.

Let $D=(V,A)$ be a directed graph with a map
$\sigma_{0}:V\to\mathbb{N}$ called a \term{chip-arrangement}. The
number $\sigma_{0}(v)$ records the \textit{number of chips} on vertex
$v$ in $\sigma_{0}$. Given a chip-arrangement~$\sigma$ a vertex $v$
can be \term{fire}\nct{d} if it contains more chips than its
out-degree, i.e.  $\sigma(v)\geq {\sf outdeg}(v)$. Firing~$v$ consist
in sending a chip along each of the out-going arcs of $v$ to their
respective end-vertices. The new chip-arrangement is called
$\sigma^{v}$. Define a directed graph $\textmd{CFG}(D,\sigma_{0})$ on
the set of chip-arrangements, which can be obtained from $\sigma_{0}$
by a firing-sequence. The arc $(\sigma,\sigma^{v})$ is naturally
colored with $v\in V$. By definition of $\textmd{CFG}(D,\sigma_{0})$
every vertex $\sigma$ lies on a directed $(\sigma_{0},\sigma)$-path.
If $\textmd{CFG}(D,\sigma_{0})$ is also acyclic the chip-firing game
is called \term{finite}. We obtain the well known:

\begin{prop}
  The states of a finite chip-firing game carry the structure of a
  ULD.
\end{prop}
\begin{proof}
  Since $\textmd{CFG}(D,\sigma_{0})$ is acyclic and every vertex lies
  on a directed $(\sigma_{0},\sigma)$-path
  $\textmd{CFG}(D,\sigma_{0})$ is connected and has a unique sink
  $\sigma_{0}$. It only remains to show, that the natural coloring of
  $\textmd{CFG}(D,\sigma_{0})$ is a U-coloring.

  Clearly we have property $\textmd{U}_1$. To prove property
  $\textmd{U}_2$ let $\sigma\in \textmd{CFG}(D,\sigma_{0})$ and
  $v,w\in V$ ready for firing in $\sigma$. Since the firing of $v$ can
  only increase the number of chips on $w$, after firing $v$ still $w$
  can be fired. This means that coloring the arcs of
  $\textmd{CFG}(D,\sigma_{0})$ by the vertices that have been fired
  leads to a U-coloring.
	
  By Theorem ~\ref{thm:main} the digraph $\textmd{CFG}(D,\sigma_{0})$
  is the cover graph of a ULD.
\end{proof}

The ULD-properties imply that in a finite game there is a unique
chip-arrangement~$\sigma^*$ where the game starting in $\sigma_0$ ends 
and that all firing sequences from a given
chip-arrangement $\sigma$ to the maximum $\sigma^*$ 
fire the same multiset of vertices.

Observe that the property which makes a U-poset a ULD, namely the
existence of a global minimum, is rather artificially achieved in the
case of chip-firing games. Instead of considering only the
arrangements that can be reached from a starting chip-arrangement
$\sigma_{0}$ by \textit{upwards} transformations, we can endow the
structure by the inverse operation of firing and call it
\textit{co-firing}.  Co-firing a vertex $v$ means sending one chip
along all the arcs $(v,w)$ from $w$ to $v$. Note that this requires
that every out-neighbour of $v$ owns a chip. We define the
\term{complete chip-firing game} $\textmd{CCFG}(D,\sigma_{0})$ to be
the digraph on all the chip-arrangements that can be reached from
$\sigma_{0}$ by firing and co-firing vertices. Under the assumption
that $\textmd{CFG}(D,\sigma_{0})$ is acyclic the obvious coloring of
the edges makes $\textmd{CCFG}(D,\sigma_{0})$ a U-poset.  This U-poset
is no lattice, but all its upsets correspond to ULDs obtained from
chip-firing games in the original sense, e.g. Figure~\ref{fig:chipfire}.

   \calc_figscale{41}
    \begin{figure}[htb]
    \centerline{\input{\path/chipfire.pstex_t}}
    \caption{\label{fig:chipfire}}
    \end{figure}
    VC
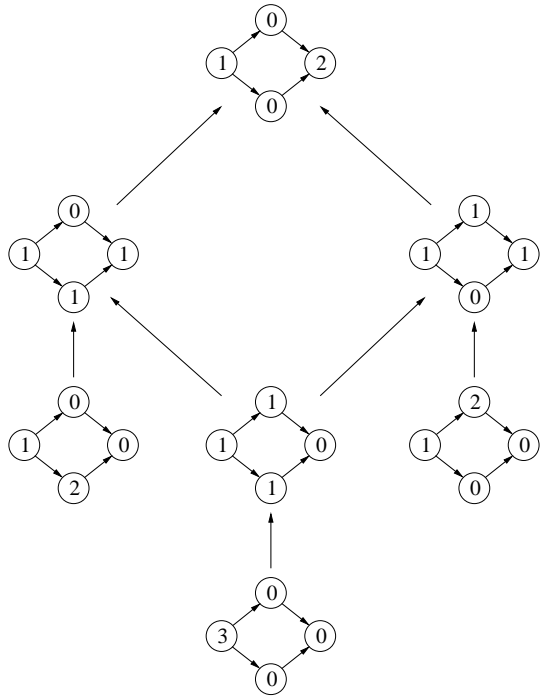
{A U-poset induced by a complete chip-firing game.}

It could be worthwhile to investigate the U-posets coming from
complete chip-firing games.
As a first step we provide the following proposition related to
Remark~2. Note that the property stated in the proposition does not
hold for the U-poset of Figure~\ref{fig:semi-uld}.

\begin{prop}
  If $(P,\leq)$ is a U-poset induced by a CCFG, then for every $s\in
  P$ there is a unique inclusion minimal set $M_{s}$ of
  meet-irreducibles such that $\bigwedge M_{s}=s$
\end{prop}

\begin{proof}
  By Remark~2. we know that for every $s\in P$ there is a unique
  inclusion minimal set $M_{s}$ of meet-irreducibles such that
  $s\in\bigwedge M_{s}$. So suppose there are distinct $s,t\in
  \bigwedge M_{s}$, i.e. $M_{s}=M_{t}$. Partition $M_{s}=U\cup V$ and let
  $u=\bigwedge U$ and $v=\bigwedge V$ and $w=u\vee v$. The relation
  between these elements is as in Figure~\ref{fig:semi-uld}.

  For two states $i < j$ in $\textmd{CCFG}(D,\sigma_{0})$ let
  $p_{i,j}$ be a directed $(i,j)$-path and $c({i,j})$ be the colorset
  of $p_{i,j}$. Recall from Section~\ref{sec:characterization of ULDs}
  that $c({i,j})$ is indeed independent of the choice of $p_{i,j}$.
  Since $s=u\wedge v$ and $w=u\vee v$ we have $c({s,u})=c({v,w})$ and
  $c({s,v})=c({u,w})$. Similarly for $t$ we obtain $c({t,u})=c({v,w})$
  and $c({t,v})=c({u,w})$.  This yields $c({s,u})=c({t,u})$ and
  $c({s,v})=c({t,v})$.  But if a chip-configuration $u\in
  \textmd{CCFG}(D,\sigma_{0})$ can be obtained from $s$ and $t$ by
  firing the same multiset of vertices then $s=t$.  This is $\bigwedge
  M_{s}=s$.
\end{proof}

\vskip-9mm
\vbox{}

\section{Concluding Remarks}
\label{sec:conclusion}

Our characterization of ULDs originates from a characterization of
matrices whose {\it flip-flop} poset generate a distributive lattice
in the diploma-thesis of Knauer~\cite{Kna-07}. It turned out that this
tool yields handy proofs for the known distributive lattices from
graphs. The extraction of the data, that were necessary for the proof,
lead to $\Delta$-bonds.  While preparing the present paper we observed
that every $\Delta$-bond lattice is isomorphic to a $\mathbf{0}$-bond
lattice, hence, the lattice structure can nicely be proved via
potentials. This lead us to the notion of a D-polytope, i.e., of a
polytope $P$ such that with points $x$ and $y$ also $\max(x,y)$ and
$\min(x,y)$ are in $P$. We can characterize the bounding hyperplanes
of D-polytoes.  This allows to associate a weighted digraph with every
D-polytope such that the {\it generalized bonds} of the weighted
digraph form a distributive lattice. This will be the topic of the
forthcoming paper~\cite{FeKn-08}.

There are two types of problems closely related to this article
where we would like to see progress:

\Bitem 
Lattices of $\Delta$-bonds depend on the choice of a
forbidden vertex $v_0 \in V$. Choosing another forbidden vertex $v_1$ 
yields a different lattice on the same set of objects. Is there an easy
description of the transformation $\PD(v_0) \to \PD(v_1)$? 

\Bitem The generation of a random element from a distributive lattice
is a nice application for {\it coupling from the past} (c.f. Propp and
Wilson~\cite{PrWi-96}). The challenge is to find good estimates for the
mixing time, see Propp~\cite{Pr-97}. What if the lattice is a 
$\Delta$-bond lattice?
\medskip

\ni
Thanks to {\'E}ric Fusy for fruitful discussions during his stay in Berlin.


\bibliography{literature,
  /homes/combi/felsner/DriMiProj/DFG/Graphor-Antrag07/antrag07-bib}
\bibliographystyle{/homes/combi/felsner/TexStuff/my-siam}


\end{document}

%% file: Figures/Urule2.pstex_t
\begin{picture}(0,0)%
\includegraphics[scale=\figscale]{Urule2.pstex}%
\end{picture}%
\setlength{\unitlength}{4144sp}%

\multiply\unitlength by \magproz
\divide\unitlength by 100

\begingroup\makeatletter\ifx\SetFigFont\undefined%
\gdef\SetFigFont#1#2#3#4#5{%
  \reset@font\fontsize{#1}{#2pt}%
  \fontfamily{#3}\fontseries{#4}\fontshape{#5}%
  \selectfont}%
\fi\endgroup%
\begin{picture}(3698,1390)(4291,-6176)
\put(4951,-6136){\makebox(0,0)[b]{\smash{{\SetFigFont{12}{14.4}{\rmdefault}{\mddefault}{\updefault}{\color[rgb]{0,0,0}$v$}%
}}}}
\put(5466,-5375){\makebox(0,0)[b]{\smash{{\SetFigFont{12}{14.4}{\rmdefault}{\mddefault}{\updefault}{\color[rgb]{0,0,0}$w$}%
}}}}
\put(4430,-5371){\makebox(0,0)[b]{\smash{{\SetFigFont{12}{14.4}{\rmdefault}{\mddefault}{\updefault}{\color[rgb]{0,0,0}$u$}%
}}}}
\put(6691,-5378){\makebox(0,0)[b]{\smash{{\SetFigFont{12}{14.4}{\rmdefault}{\mddefault}{\updefault}{\color[rgb]{0,0,0}$u$}%
}}}}
\put(7828,-5382){\makebox(0,0)[b]{\smash{{\SetFigFont{12}{14.4}{\rmdefault}{\mddefault}{\updefault}{\color[rgb]{0,0,0}$w$}%
}}}}
\put(7291,-6110){\makebox(0,0)[b]{\smash{{\SetFigFont{12}{14.4}{\rmdefault}{\mddefault}{\updefault}{\color[rgb]{0,0,0}$v$}%
}}}}
\put(7291,-4894){\makebox(0,0)[b]{\smash{{\SetFigFont{12}{14.4}{\rmdefault}{\mddefault}{\updefault}{\color[rgb]{0,0,0}$z$}%
}}}}
\end{picture}%

%% file: Figures/path-chasing.pstex_t
\begin{picture}(0,0)%
\includegraphics[scale=\figscale]{path-chasing.pstex}%
\end{picture}%
\setlength{\unitlength}{4144sp}%

\multiply\unitlength by \magproz
\divide\unitlength by 100

\begingroup\makeatletter\ifx\SetFigFont\undefined%
\gdef\SetFigFont#1#2#3#4#5{%
  \reset@font\fontsize{#1}{#2pt}%
  \fontfamily{#3}\fontseries{#4}\fontshape{#5}%
  \selectfont}%
\fi\endgroup%
\begin{picture}(8172,3189)(-2518,-9466)
\put(2466,-9421){\makebox(0,0)[lb]{\smash{{\SetFigFont{12}{14.4}{\rmdefault}{\mddefault}{\updefault}{\color[rgb]{0,0,0}$x_0=x$}%
}}}}
\put(2896,-9031){\makebox(0,0)[lb]{\smash{{\SetFigFont{12}{14.4}{\rmdefault}{\mddefault}{\updefault}{\color[rgb]{0,0,0}$x_1$}%
}}}}
\put(3626,-8221){\makebox(0,0)[lb]{\smash{{\SetFigFont{12}{14.4}{\rmdefault}{\mddefault}{\updefault}{\color[rgb]{0,0,0}$x_{\ell-1}$}%
}}}}
\put(3986,-7871){\makebox(0,0)[lb]{\smash{{\SetFigFont{12}{14.4}{\rmdefault}{\mddefault}{\updefault}{\color[rgb]{0,0,0}$x_\ell$}%
}}}}
\put(4296,-7561){\makebox(0,0)[lb]{\smash{{\SetFigFont{12}{14.4}{\rmdefault}{\mddefault}{\updefault}{\color[rgb]{0,0,0}$x_{\ell+1}=y_\ell$}%
}}}}
\put(4916,-6881){\makebox(0,0)[lb]{\smash{{\SetFigFont{12}{14.4}{\rmdefault}{\mddefault}{\updefault}{\color[rgb]{0,0,0}$x_k = z$}%
}}}}
\put(2816,-7647){\makebox(0,0)[rb]{\smash{{\SetFigFont{12}{14.4}{\rmdefault}{\mddefault}{\updefault}{\color[rgb]{0,0,0}$y_{\ell-1}$}%
}}}}
\put(2046,-8337){\makebox(0,0)[rb]{\smash{{\SetFigFont{12}{14.4}{\rmdefault}{\mddefault}{\updefault}{\color[rgb]{0,0,0}$y_1$}%
}}}}
\put(1586,-8867){\makebox(0,0)[rb]{\smash{{\SetFigFont{12}{14.4}{\rmdefault}{\mddefault}{\updefault}{\color[rgb]{0,0,0}$y=y_0$}%
}}}}
\put(1473,-6881){\makebox(0,0)[lb]{\smash{{\SetFigFont{12}{14.4}{\rmdefault}{\mddefault}{\updefault}{\color[rgb]{0,0,0}$x_k = z$}%
}}}}
\put(-547,-9031){\makebox(0,0)[lb]{\smash{{\SetFigFont{12}{14.4}{\rmdefault}{\mddefault}{\updefault}{\color[rgb]{0,0,0}$x_1$}%
}}}}
\put(-977,-9421){\makebox(0,0)[lb]{\smash{{\SetFigFont{12}{14.4}{\rmdefault}{\mddefault}{\updefault}{\color[rgb]{0,0,0}$x_0=x$}%
}}}}
\put(-1857,-8867){\makebox(0,0)[rb]{\smash{{\SetFigFont{12}{14.4}{\rmdefault}{\mddefault}{\updefault}{\color[rgb]{0,0,0}$y=y_0$}%
}}}}
\put(-1397,-8337){\makebox(0,0)[rb]{\smash{{\SetFigFont{12}{14.4}{\rmdefault}{\mddefault}{\updefault}{\color[rgb]{0,0,0}$y_1$}%
}}}}
\put(651,-6397){\makebox(0,0)[rb]{\smash{{\SetFigFont{12}{14.4}{\rmdefault}{\mddefault}{\updefault}{\color[rgb]{0,0,0}$y_k$}%
}}}}
\put(-949,-6807){\makebox(0,0)[b]{\smash{{\SetFigFont{12}{14.4}{\rmdefault}{\mddefault}{\updefault}{\color[rgb]{0,0,0}(a)}%
}}}}
\put(3361,-6807){\makebox(0,0)[b]{\smash{{\SetFigFont{12}{14.4}{\rmdefault}{\mddefault}{\updefault}{\color[rgb]{0,0,0}(b)}%
}}}}
\end{picture}%

%% file: Figures/counterE.pstex_t
\begin{picture}(0,0)%
\includegraphics[scale=\figscale]{counterE.pstex}%
\end{picture}%
\setlength{\unitlength}{4144sp}%

\multiply\unitlength by \magproz
\divide\unitlength by 100

\begingroup\makeatletter\ifx\SetFigFont\undefined%
\gdef\SetFigFont#1#2#3#4#5{%
  \reset@font\fontsize{#1}{#2pt}%
  \fontfamily{#3}\fontseries{#4}\fontshape{#5}%
  \selectfont}%
\fi\endgroup%
\begin{picture}(4760,3685)(-42,-13897)
\put(711,-10710){\makebox(0,0)[rb]{\smash{{\SetFigFont{12}{14.4}{\rmdefault}{\mddefault}{\updefault}{\color[rgb]{0,0,0}$y = x_2$}%
}}}}
\put(2498,-10860){\makebox(0,0)[b]{\smash{{\SetFigFont{12}{14.4}{\rmdefault}{\mddefault}{\updefault}{\color[rgb]{0,0,0}$x_3$}%
}}}}
\put(4257,-10760){\makebox(0,0)[lb]{\smash{{\SetFigFont{12}{14.4}{\rmdefault}{\mddefault}{\updefault}{\color[rgb]{0,0,0}$x_4$}%
}}}}
\put(3284,-12200){\makebox(0,0)[lb]{\smash{{\SetFigFont{12}{14.4}{\rmdefault}{\mddefault}{\updefault}{\color[rgb]{0,0,0}$x_5$}%
}}}}
\put(2112,-11961){\makebox(0,0)[lb]{\smash{{\SetFigFont{12}{14.4}{\rmdefault}{\mddefault}{\updefault}{\color[rgb]{0,0,0}$x_1$}%
}}}}
\put(2470,-13852){\makebox(0,0)[b]{\smash{{\SetFigFont{12}{14.4}{\rmdefault}{\mddefault}{\updefault}{\color[rgb]{0,0,0}$x_0=x$}%
}}}}
\end{picture}%

%% file: Figures/semi-uld.pstex_t
\begin{picture}(0,0)%
\includegraphics[scale=\figscale]{semi-uld.pstex}%
\end{picture}%
\setlength{\unitlength}{4144sp}%

\multiply\unitlength by \magproz
\divide\unitlength by 100

\begingroup\makeatletter\ifx\SetFigFont\undefined%
\gdef\SetFigFont#1#2#3#4#5{%
  \reset@font\fontsize{#1}{#2pt}%
  \fontfamily{#3}\fontseries{#4}\fontshape{#5}%
  \selectfont}%
\fi\endgroup%
\begin{picture}(1813,2083)(1794,-11981)
\put(2071,-10951){\makebox(0,0)[rb]{\smash{{\SetFigFont{12}{14.4}{\rmdefault}{\mddefault}{\updefault}{\color[rgb]{0,0,0}$u$}%
}}}}
\put(3331,-10951){\makebox(0,0)[lb]{\smash{{\SetFigFont{12}{14.4}{\rmdefault}{\mddefault}{\updefault}{\color[rgb]{0,0,0}$v$}%
}}}}
\put(3331,-11941){\makebox(0,0)[b]{\smash{{\SetFigFont{12}{14.4}{\rmdefault}{\mddefault}{\updefault}{\color[rgb]{0,0,0}$t$}%
}}}}
\put(2071,-11941){\makebox(0,0)[b]{\smash{{\SetFigFont{12}{14.4}{\rmdefault}{\mddefault}{\updefault}{\color[rgb]{0,0,0}$s$}%
}}}}
\put(2701,-10006){\makebox(0,0)[b]{\smash{{\SetFigFont{12}{14.4}{\rmdefault}{\mddefault}{\updefault}{\color[rgb]{0,0,0}$w$}%
}}}}
\end{picture}%

%% file: Figures/chipfire.pstex_t
\begin{picture}(0,0)%
\includegraphics[scale=\figscale]{chipfire.pstex}%
\end{picture}%
\setlength{\unitlength}{4144sp}%

\multiply\unitlength by \magproz
\divide\unitlength by 100

\begingroup\makeatletter\ifx\SetFigFont\undefined%
\gdef\SetFigFont#1#2#3#4#5{%
  \reset@font\fontsize{#1}{#2pt}%
  \fontfamily{#3}\fontseries{#4}\fontshape{#5}%
  \selectfont}%
\fi\endgroup%
\begin{picture}(7756,10096)(1568,-8484)
\end{picture}%